# THE EFFICIENCY OF THE ESTIMATORS OF THE PARAMETERS IN GARCH PROCESSES

By István Berkes[1] and Lajos Horváth[2]

*Hungarian Academy of Sciences and University of Utah*

We propose a class of estimators for the parameters of a GARCH$(p,q)$ sequence. We show that our estimators are consistent and asymptotically normal under mild conditions. The quasi-maximum likelihood and the likelihood estimators are discussed in detail. We show that the maximum likelihood estimator is optimal. If the tail of the distribution of the innovations is polynomial, even a quasi-maximum likelihood estimator based on exponential density performs better than the standard normal density-based quasi-likelihood estimator of Lee and Hansen and Lumsdaine.

**1. Introduction.** The generalized autoregressive conditional heteroscedastic (GARCH) process was introduced by Bollerslev (1986). The GARCH process has received considerable attention from applied as well as from theoretical points of view. We say that $\{y_k, -\infty < k < \infty\}$ is a GARCH$(p,q)$ process if it satisfies the equations

$$y_k = \sigma_k \varepsilon_k \tag{1.1}$$

and

$$\sigma_k^2 = \omega + \sum_{1 \leq i \leq p} \alpha_i y_{k-i}^2 + \sum_{1 \leq j \leq q} \beta_j \sigma_{k-j}^2, \tag{1.2}$$

where

$$\omega > 0, \ \alpha_i \geq 0, \ 1 \leq i \leq p, \ \beta_j \geq 0, \ 1 \leq j \leq q \tag{1.3}$$

Received November 2001; revised October 2002.
[1]Supported by the Hungarian National Foundation for Scientific Research Grant T 29621 and NSF Grant INT-02-23262.
[2]Supported by the DOMUS Hungarica Sci. Art and NSF Grant INT-02-23262.

are constants. We also assume that

(1.4) $\{\varepsilon_i, -\infty < i < \infty\}$ are independent, identically distributed random variables.

Throughout this paper we assume that (1.1)–(1.4) hold.

The GARCH$(1,1)$ model was studied by Nelson (1991) who showed that (1.1) and (1.2) have a unique stationary solution if and only if $E\log(\beta_1 + \alpha_1\varepsilon_0^2) < 0$. The general case was investigated by Bougerol and Picard (1992a, b). Let

$$\boldsymbol{\tau}_n = (\beta_1 + \alpha_1\varepsilon_n^2, \beta_2, \ldots, \beta_{q-1}) \in \mathbb{R}^{q-1},$$
$$\boldsymbol{\xi}_n = (\varepsilon_n^2, 0, \ldots, 0) \in \mathbb{R}^{q-1}$$

and

$$\boldsymbol{\alpha} = (\alpha_2, \ldots, \alpha_{p-1}) \in \mathbb{R}^{p-2}.$$

[Clearly, without loss of generality, we may and shall assume $\min(p,q) \geq 2$.] Define the $(p+q-1) \times (p+q-1)$ matrix $A_n$, written in block form, by

$$A_n = \begin{bmatrix} \boldsymbol{\tau}_n & \beta_q & \boldsymbol{\alpha} & \alpha_p \\ \mathbf{I}_{q-1} & 0 & 0 & 0 \\ \boldsymbol{\xi}_n & 0 & 0 & 0 \\ 0 & 0 & \mathbf{I}_{p-2} & 0 \end{bmatrix},$$

where $\mathbf{I}_{q-1}$ and $\mathbf{I}_{p-2}$ are the identity matrices of size $q-1$ and $p-2$, respectively. The norm of any $d \times d$ matrix $M$ is defined by

$$\|M\| = \sup\{\|M\mathbf{x}\|_d/\|\mathbf{x}\|_d : \mathbf{x} \in \mathbb{R}^d, \mathbf{x} \neq \mathbf{0}\},$$

where $\|\cdot\|_d$ is the usual (Euclidean) norm in $\mathbb{R}^d$. The top Liapounov exponent $\gamma_L$ associated with the sequence $\{A_n, -\infty < n < \infty\}$ is

$$\gamma_L = \inf_{1 \leq n < \infty} \frac{1}{n+1} E\log\|A_0 A_1 \cdots A_n\|,$$

assuming that

(1.5) $$E(\log\|A_0\|) < \infty.$$

Bougerol and Picard (1992a, b) showed that if (1.5) holds, then (1.1) and (1.2) have a unique stationary solution if and only if

(1.6) $$\gamma_L < 0.$$

The estimation of the parameter $\boldsymbol{\theta} = (\omega, \alpha_1, \ldots, \alpha_p, \beta_1, \ldots, \beta_q)$ has been studied by several authors. Lee and Hansen (1994) and Lumsdaine (1996) used the quasi-maximum likelihood method to estimate the parameters from the sample $y_1, \ldots, y_n$ in GARCH$(1,1)$ models. The idea behind the quasi-maximum likelihood method is the following. The likelihood function is derived under the assumption that $\varepsilon_0$ is standard normal. The estimator is the point where the likelihood function reaches its largest value. The estimator in Lee and Hansen (1994) and Lumsdaine (1996) is "local" since the likelihood function is maximized in a small neighborhood of $\boldsymbol{\theta}$. They show that the quasi-maximum likelihood estimator is consistent and asymptotically normal without assuming the normality of $\varepsilon_0$. However, very strict conditions are assumed on the distribution of $\varepsilon_0$ and the value of $\boldsymbol{\theta}$. Berkes, Horváth and Kokoszka (2003) investigated the asymptotic properties of the quasi-maximum likelihood estimator for $\boldsymbol{\theta}$ in GARCH$(p,q)$ models. Berkes, Horváth and Kokoszka (2003) obtained their asymptotic results under weak conditions. Berkes and Horváth (2003) showed that the quasi-maximum likelihood estimator cannot be $n^{-1/2}$-consistent if $E|\varepsilon_0|^\kappa = \infty$ for some $0 < \kappa < 4$. This shows the limitations of the quasi-maximum likelihood estimation method. The existence of the GARCH$(p,q)$ sequence requires only that $E|\log \varepsilon_0^2| < \infty$ but the estimation works only if $E|\varepsilon_0|^\kappa < \infty$ with some $\kappa > 4$. The quasi-maximum likelihood estimator does not use the distribution of $\varepsilon_0$ and therefore, as we shall see, it is not efficient. If $E\varepsilon_0 = 0$ and $E\varepsilon_0^2 = 1$, then $\sigma_k^2$ is the conditional variance of $y_k$ given the past. However, without any moment conditions, $\sigma_k$ is the conditional scaling parameter of $y_k$.

Since $\sigma_k$ is defined by a recursion, we use a recursion to define our estimator. Let $\mathbf{u} = (x, \mathbf{s}, \mathbf{t}) \in \mathbb{R}^{p+q+1}$, $x \in \mathbb{R}$, $\mathbf{s} \in \mathbb{R}^p$ and $\mathbf{t} \in \mathbb{R}^q$. We start with the initial conditions: if $q \geq p$, then

$$c_0(\mathbf{u}) = x/(1 - (t_1 + \cdots + t_q)),$$
$$c_1(\mathbf{u}) = s_1,$$
$$c_2(\mathbf{u}) = s_2 + t_1 c_1(\mathbf{u}),$$
$$\vdots$$
$$c_p(\mathbf{u}) = s_p + t_1 c_{p-1}(\mathbf{u}) + \cdots + t_{p-1} c_1(\mathbf{u}),$$
$$c_{p+1}(\mathbf{u}) = t_1 c_p(\mathbf{u}) + \cdots + t_p c_1(\mathbf{u}),$$
$$\vdots$$
$$c_q(\mathbf{u}) = t_1 c_{q-1}(\mathbf{u}) + \cdots + t_{q-1} c_1(\mathbf{u}),$$

and if $q < p$, the equations above are replaced with

$$c_0(\mathbf{u}) = x/(1 - (t_1 + \cdots + t_q)),$$
$$c_1(\mathbf{u}) = s_1,$$
$$c_2(\mathbf{u}) = s_2 + t_1 c_1(\mathbf{u}),$$

$$\vdots$$
$$c_{q+1}(\mathbf{u}) = s_{q+1} + t_1 c_q(\mathbf{u}) + \cdots + t_q c_1(\mathbf{u}),$$
$$\vdots$$
$$c_p(\mathbf{u}) = s_p + t_1 c_{p-1}(\mathbf{u}) + \cdots + t_q c_{p-q}(\mathbf{u}).$$

In general, if $i > R = \max(p,q)$, then

(1.7) $$c_i(\mathbf{u}) = t_1 c_{i-1}(\mathbf{u}) + t_2 c_{i-2}(\mathbf{u}) + \cdots + t_q c_{i-q}(\mathbf{u}).$$

We choose an arbitrary positive function $h$ and define

$$\hat{L}_n(\mathbf{u}) = \frac{1}{n} \sum_{1 < k \leq n} \log\left\{ \frac{1}{\hat{w}_k^{1/2}(\mathbf{u})} h(y_k/\hat{w}_k^{1/2}(\mathbf{u})) \right\},$$

where

$$\hat{w}_k(\mathbf{u}) = c_0(\mathbf{u}) + \sum_{1 \leq i < k} c_i(\mathbf{u}) y_{k-i}^2.$$

Let $0 < \underline{u} < \overline{u}$, $0 < \rho_0 < 1$, $q\underline{u} < \rho_0$ and define

$$U = \{\mathbf{u} : t_1 + t_2 + \cdots + t_q \leq \rho_0 \text{ and}$$
$$\underline{u} \leq \min(x, s_1, s_2, \ldots, s_p, t_1, t_2, \ldots, t_q)$$
$$\leq \max(x, s_1, s_2, \ldots, s_p, t_1, t_2, \ldots, t_q) \leq \overline{u}\}.$$

From now on we replace (1.3) with the somewhat stronger condition

(1.8) $\boldsymbol{\theta}$ is in the interior of $U$.

We use $|\cdot|$ to denote the maximum norm of vectors and matrices. Let $x \vee y = \max(x, y)$. In this paper we study the asymptotic properties of

$$\hat{\boldsymbol{\theta}}_n = \arg\max_{\mathbf{u} \in U} \hat{L}_n(\mathbf{u}).$$

We note that $\hat{L}_n(\mathbf{u})$ is a continuously differentiable function, so standard numerical methods can be used to compute $\hat{\boldsymbol{\theta}}_n$.

In our first result we give a sufficient criterion for $|\hat{\boldsymbol{\theta}}_n - \boldsymbol{\theta}| \to 0$ a.s. To state this result we will need some additional regularity conditions:

(1.9) the polynomials $\alpha_1 x + \alpha_2 x^2 + \cdots + \alpha_p x^p$ and $1 - \beta_1 x - \beta_2 x^2 - \cdots - \beta_q x^q$ are coprimes in the set of polynomials with real coefficients,

(1.10) $\varepsilon_0^2$ is a nondegenerate random variable

and

(1.11) $$\lim_{t \to 0} t^{-\mu} P\{\varepsilon_0^2 \le t\} = 0, \qquad \text{with some } \mu > 0.$$

Condition (1.8) is somewhat stronger than (1.3) but $\beta_1 + \cdots + \beta_q < 1$ is a necessary condition for the existence of a GARCH$(p, q)$ sequence [cf. Berkes, Horváth and Kokoszka (2003)]. Assumptions (1.9) and (1.10) are needed to uniquely identify the parameter $\boldsymbol{\theta}$. So far all our conditions are related to the structure of the GARCH$(p, q)$ process. The following set of conditions concerns the moments of $\varepsilon_0$ and the smoothness of $h$:

(1.12) $$E|\varepsilon_0^2|^\kappa < \infty \qquad \text{with some } \kappa > 0,$$

and there is $0 < C_0 < \infty$ such that

(1.13) $$E|\log h(\varepsilon_0 t)| \le C_0(t^{\nu_0} + 1) \qquad \text{for all } t > 0, \text{ with some } 0 \le \nu_0 < 2\kappa.$$

Let

$$g(y, t) = \log\{th(yt)\}, \qquad -\infty < y < \infty, \; t > 0,$$

and

$$g_1(y, t) = \frac{\partial}{\partial t} g(y, t), \qquad -\infty < y < \infty, \; t > 0.$$

We also assume that there is a function $C_1(y)$ such that

(1.14) $$|g_1(y, t)| \le C_1(y)(t^{\nu_1} + 1)/t \qquad \text{for all } 0 < t < \infty \text{ and } y \in \mathbb{R},$$
$$\text{with some } 0 \le \nu_1 < 2\kappa,$$

and

(1.15) $$EC_1(\varepsilon_0) < \infty.$$

If $h$ is a density, then condition (1.14) means that the density function $th(yt)$ is smooth in the parameter $t$.

We will show in Lemma 4.1 that

$$L(\mathbf{u}) = E \log\left\{ \frac{1}{(w_0(\mathbf{u}))^{1/2}} h(y_0/(w_0(\mathbf{u}))^{1/2}) \right\}$$

exists for all $\mathbf{u} \in U$, where

$$w_k(\mathbf{u}) = c_0(\mathbf{u}) + \sum_{1 \le i < \infty} c_i(\mathbf{u}) y_{k-i}^2.$$

We note that

$$w_k(\boldsymbol{\theta}) = \sigma_k^2.$$

The following condition will imply [see (4.6)] that $L(\mathbf{u})$ has a unique maximum in $U$ at $\boldsymbol{\theta}$:

(1.16) $$Eg(\varepsilon_0, t) < Eg(\varepsilon_0, 1) \qquad \text{for all } 0 < t < \infty, \; t \ne 1.$$

THEOREM 1.1. *If (1.5), (1.6) and (1.8)–(1.16) hold, then*

$$\hat{\boldsymbol{\theta}}_n \to \boldsymbol{\theta} \qquad a.s.$$

The proof of Theorem 1.1 will be given in Section 4.

Next we discuss the asymptotic normality of $n^{1/2}(\hat{\boldsymbol{\theta}}_n - \boldsymbol{\theta})$. We need further smoothness conditions on $th(yt)$. Let $g_2(y,t)$ and $g_3(y,t)$ be the second and third derivatives of $g(y,t)$ with respect to $t$. We assume that there are functions $C_2$ and $C_3$ such that

$$(1.17) \quad |g_2(y,t)| \leq C_2(y)(t^{\nu_2} + 1)/t^2 \qquad \text{for all } 0 < t < \infty \text{ and } y \in \mathbb{R},$$
$$\text{with some } 0 \leq \nu_2 < \infty,$$

$$(1.18) \quad EC_2(\varepsilon_0) < \infty,$$

$$(1.19) \quad |g_3(y,t)| \leq C_3(y)(t^{\nu_3} + 1)/t^3 \qquad \text{for all } 0 < t < \infty \text{ and } y \in \mathbb{R},$$
$$\text{with some } 0 \leq \nu_3 < \infty,$$

and

$$(1.20) \quad EC_3(\varepsilon_0) < \infty.$$

We use $w_k'(\mathbf{u})$ to denote the row vector of the derivatives of $w_k(\mathbf{u})$ and $w_k''(\mathbf{u})$ the matrix of the second-order partial derivatives of $w_k(\mathbf{u})$ (the Hessian matrix). Berkes, Horváth and Kokoszka (2003) showed that

$$\mathbf{A} = E(w_0'(\boldsymbol{\theta})/w_0(\boldsymbol{\theta}))^T (w_0'(\boldsymbol{\theta})/w_0(\boldsymbol{\theta}))$$

exists and is nonsingular ($T$ denotes the transpose). We also assume that

$$(1.21) \quad 0 < Eg_1^2(\varepsilon_0, 1) < \infty,$$

$$(1.22) \quad E|g_2(\varepsilon_0, 1)| < \infty \quad \text{and} \quad Eg_2(\varepsilon_0, 1) \neq 0.$$

If (1.21) and (1.22) hold, then

$$0 < \tau^2 = \frac{Eg_1^2(\varepsilon_0, 1)}{(Eg_2(\varepsilon_0, 1))^2} < \infty.$$

The multivariate normal distribution with mean $\mathbf{0}$ and covariance matrix $\mathbf{D}$ will be denoted by $\mathbf{N}(\mathbf{0}, \mathbf{D})$.

THEOREM 1.2. *If (1.5), (1.6) and (1.8)–(1.22) hold, then*

$$n^{1/2}(\hat{\boldsymbol{\theta}}_n - \boldsymbol{\theta}) \xrightarrow{\mathcal{D}} \mathbf{N}(\mathbf{0}, 4\tau^2 \mathbf{A}^{-1}).$$

This result will be proven in Section 4.

REMARK 1.1. Let $f(y)$ denote the density function of $\varepsilon_0$ and $I_f(t)$ the Fisher information number of the scale family $tf(xt)$, $t > 0$. If $\varepsilon_1, \ldots, \varepsilon_n$ is known, then $\hat{t}_n = \arg\max\{\prod_{1 \leq i \leq n} tf(\varepsilon_i t) : t > 0\}$ can be used to estimate the scale parameter. One can verify that under suitable regularity conditions $n^{1/2}(\hat{t}_n - 1)$ will be asymptotically normal with mean 0 and variance $\tau^2$. So by Lehmann [(1991), page 406] we conclude that

$$\tau^2 \geq \frac{1}{I_f(1)}, \tag{1.23}$$

and we have the equality in (1.23) when $h = f$.

REMARK 1.2. Newey and Steigerwald (1997) consider more general models which include the GARCH$(p, q)$ sequence. They point out that identification of the parameters in the drift term might be difficult. In our paper we study the estimation of the parameters in the error process of the Newey–Steigerwald model. This is the part which makes GARCH different from other time series. Our results cannot be applied directly to other versions of GARCH but our method can be used to investigate the properties of estimators in LGARCH [Bollerslev (1986)], NGARCH [Engle and Ng (1993)], MGARCH [Geweke (1986)], EGARCH [Nelson (1991)] and VGARCH [Engle and Ng (1993)].

REMARK 1.3. Lee and Hansen (1994) assume that the observed sequence $y_k$ is a stationary and ergodic martingale. They also assume that

$$Ey_0^2 < \infty.$$

We do not impose this moment condition. Under our conditions we have only that

$$E|y_0|^\delta < \infty, \quad \text{with some } \delta > 0.$$

It would be interesting and practically useful to extend the results of Lee and Hansen (1994) to the present situation.

REMARK 1.4. Drost and Klaassen (1997) showed that there is a reparametrization of GARCH$(1, 1)$ such that the efficient score functions in the parametric model of the autoregression parameters are orthogonal to the tangent space generated by the nuisance parameter, thus suggesting that adaptive estimation of the parameters is possible. Drost and Klaassen (1997) construct adaptive and hence efficient estimators in the reparametrized GARCH$(1, 1)$ in a mean-type context.

Next we consider three special choices of $h$.

## 2. Examples.

EXAMPLE 2.1. Let $h(t) = (2\pi)^{-1/2} \exp(-t^2/2)$ (the standard normal density function). Using this function in the definition of $\hat{L}_n$, we get the quasi-maximum likelihood estimator investigated by Lee and Hansen (1994) and Lumsdaine (1996). Elementary calculations give that $|\log h(yt)| \leq C_0(y^2t^2 + 1)$,

$$g_1(y,t) = (1 - y^2t^2)/t, \qquad |g_1(t)| \leq (1+y^2)(1+t^2)/t,$$
$$g_2(y,t) = -(1 + y^2t^2)/t^2, \qquad |g_2(t)| \leq (1+y^2)(1+t^2)/t^2$$

and $g_3(y,t) = 2/t^3$. It is easy to see that $t = 1/E\varepsilon_0^2$ is the unique solution of the equation $Eg_1(\varepsilon_0, t) = 0$ and $Eg(\varepsilon_0, t)$ has a unique maximum at $1/E\varepsilon_0^2$. If we assume that

$$(2.1) \qquad E\varepsilon_0^2 = 1,$$

then condition (1.16) is satisfied. We note that (2.1) is a standard condition assumed by Lee and Hansen (1994) and Lumsdaine (1996). Clearly, $g_1(\varepsilon_0, 1) = 1 - \varepsilon_0^2$ and $g_2(\varepsilon_0, 1) = -1 - \varepsilon_0^2$. Hence (1.21) holds if and only if $E\varepsilon_0^4 < \infty$. Also, $Eg_2(\varepsilon_0, 1) = -2$ by (2.1) and $\tau^2 = E(1 - \varepsilon_0^2)^2/4 = (E\varepsilon_0^4 - 1)/4$. Hence the quasi-maximum likelihood estimator is almost sure consistent if $E|\varepsilon_0^2|^\kappa < \infty$ with some $\kappa > 1$ and asymptotically normal if $E\varepsilon_0^4 < \infty$.

EXAMPLE 2.2. Let $h(t) = (1/2) \exp(-|t|)$ (two-sided exponential distribution). Elementary calculations show that $|\log h(yt)| \leq 1 + |y|t$, $E|\log h(\varepsilon_0 t)| \leq 1 + tE|\varepsilon_0|$, $g(y,t) = \log t - \log 2 - |y|t$,

$$g_1(y,t) = (1 - |y|t)/t, \qquad |g_1(t)| \leq (1 + |y|)(1+t)/t,$$

$g_2(y,t) = -1/t^2$ and $g_3(y,t) = 2/t^3$. Hence the unique solution of the equation $Eg_1(\varepsilon_0, t) = 0$ is $t = 1/E|\varepsilon_0|$, which will be 1 if and only if $E|\varepsilon_0| = 1$. Assuming that

$$(2.2) \qquad E|\varepsilon_0| = 1,$$

we get that (1.16) holds. Clearly, $g_1(\varepsilon_0, 1) = 1 - |\varepsilon_0|$ and $g_2(\varepsilon_0, 1) = -1$. Hence (1.22) is always satisfied and (1.21) holds if and only if $E\varepsilon_0^2 < \infty$. Also, $\tau^2 = E(1 - |\varepsilon_0|)^2 = E\varepsilon_0^2 - 1$. Hence the exponential density based estimator is almost sure consistent if $E|\varepsilon_0^2|^\kappa < \infty$ with some $\kappa > 1/2$ and asymptotically normal if $E\varepsilon_0^2 < \infty$.

EXAMPLE 2.3. Let $h(t) = \{(\vartheta - 1)/2\}(1 + |t|)^{-\vartheta}$ with some $\vartheta > 1$. We note that $E|\log h(\varepsilon_0 t)| \leq C_0(E\log(|\varepsilon_0| + 1) + \log(1 + t) + 1)$,

$$g(y,t) = \log t + \log((\vartheta - 1)/2) - \vartheta \log(1 + |y|t),$$
$$g_1(y,t) = 1/t - \vartheta|y|/(1 + |y|t), \qquad |g_1(y,t)| \leq C_1/t,$$
$$g_2(y,t) = -\frac{1}{t^2} + \vartheta t \frac{y^2}{(1 + |y|t)^2}, \qquad |g_2(y,t)| \leq C_2/t^2$$

and $|g_3(y,t)| \leq C_3/t^3$. The unique solution of the equation $Eg_1(\varepsilon_0, t) = 0$ is $t = 1$ if and only if

$$E\left(\frac{|\varepsilon_0|}{1+|\varepsilon_0|}\right) = \frac{1}{\vartheta} \quad (2.3)$$

and since $g_2(y,t) < 0$, $Eg(\varepsilon_0, t)$ has a unique maximum at $t = 1$; that is, (1.16) holds. By (1.10) we have (1.21) and

$$Eg_2(\varepsilon_0, 1) = -1 + \vartheta E\left(\frac{|\varepsilon_0|}{1+|\varepsilon_0|}\right)^2 < -1 + \vartheta E\left(\frac{|\varepsilon_0|}{1+|\varepsilon_0|}\right) = 0,$$

showing that (1.22) holds. Thus we can estimate $\boldsymbol{\theta}$ using this $h$ as long as $E|\varepsilon_0^2|^\kappa < \infty$ with some $\kappa > 0$.

EXAMPLE 2.4. Let $h(t) = f(t)$, where $f$ is the density function of $\varepsilon_0$. Since $-\log$ is strictly convex, Jensen's inequality shows that

$$E\log\{tf(\varepsilon_0 t)/f(\varepsilon_0)\} < \log E\{tf(\varepsilon_0 t)/f(\varepsilon_0)\} = 0 \quad (2.4)$$

if

$$tf(\varepsilon_0 t)/f(\varepsilon_0) \quad \text{is nonconstant.} \quad (2.5)$$

If, following Lehmann [(1991), page 409], we assume that the distributions determined by the scale family of densities $tf(yt)$, $t > 0$, are distinct, then (2.5) holds, with the exception of $t = 1$, and therefore (1.16) holds. Also, $g_1(\varepsilon_0, 1) = 1 + \varepsilon_0 f'(\varepsilon_0)/f(\varepsilon_0)$,

$$g_2(\varepsilon_0, 1) = -1 + \varepsilon_0^2 \left\{ \frac{f''(\varepsilon_0)}{f(\varepsilon_0)} - \left(\frac{f'(\varepsilon_0)}{f(\varepsilon_0)}\right)^2 \right\}$$

and $\tau^2 = 1/I_f(1)$, where $I_f(t)$ is the Fisher information number of the scale family $tf(yt)$, $t > 0$. In this case (1.13)–(1.19) are analogous to the conditions used by Lehmann (1991), Section 6.2, to establish the asymptotic normality of the maximum likelihood estimator of the scale parameter of the family $tf(yt)$ based on independent, identically distributed observations.

Condition (1.16) connects $h$ and the distribution of the innovations. We have seen in Example 2.4 that (1.16) is always satisfied if the maximum likelihood method is used. However, using another $h$, we may have to scale the model [cf. (2.1)–(2.3)]. Next we study the effect of scaling on the estimators and their asymptotic distributions. Let us assume that our model is

$$y_k = \tilde{\sigma}_k \tilde{\varepsilon}_k, \quad (2.6)$$

$$\tilde{\sigma}_k^2 = \tilde{\omega} + \sum_{1 \leq i \leq p} \tilde{\alpha}_i y_{k-i}^2 + \sum_{1 \leq j \leq q} \tilde{\beta}_j \tilde{\sigma}_{k-j}^2. \quad (2.7)$$

The parameter of (2.6) and (2.7) is $\tilde{\boldsymbol{\theta}} = (\tilde{\omega}, \tilde{\alpha}_1, \ldots, \tilde{\alpha}_p, \tilde{\beta}_1, \ldots, \tilde{\beta}_q)$. The scaling of $\tilde{\varepsilon}_k$ will result in $\varepsilon_k = \tilde{\varepsilon}_k/d$, $d > 0$ and $\sigma_k = d\tilde{\sigma}_k$. Thus (1.1) and (1.2) hold with $\boldsymbol{\theta} = (d^2\tilde{\omega}, d^2\tilde{\alpha}_1, \ldots, d^2\tilde{\alpha}_p, \tilde{\beta}_1, \ldots, \tilde{\beta}_q)$. We choose $d$ such that (1.16) holds. By Theorem 2.2 we have that

$$n^{1/2}(\hat{\boldsymbol{\theta}}_n - \boldsymbol{\theta}) \xrightarrow{\mathcal{D}} \mathbf{N}(\mathbf{0}, 4\tau^2 \mathbf{A}^{-1}).$$

The definitions of $c_i(\mathbf{u}), 0 \leq i < \infty$, yield that

$$\frac{w'_k(\boldsymbol{\theta})}{w_k(\boldsymbol{\theta})} = \mathbf{M} \frac{w'_k(\tilde{\boldsymbol{\theta}})}{w_k(\tilde{\boldsymbol{\theta}})},$$

where $\mathbf{M} = \{M(i,j), 0 \leq i, j \leq p+q\}$, $M(i,j) = 0$ if $i \neq j$, $M(i,i) = 1/d^2$ if $0 \leq i \leq p$ and $M(i,i) = 1$ if $p < i \leq p+q$. Hence

$$n^{1/2}(\hat{\boldsymbol{\theta}}_n - \boldsymbol{\theta}) \xrightarrow{\mathcal{D}} \mathbf{N}\left(\mathbf{0}, 4\tau^2 \mathbf{M}^{-1} E\left(\frac{w'_k(\tilde{\boldsymbol{\theta}})}{w_k(\tilde{\boldsymbol{\theta}})}\right)^T \left(\frac{w'_k(\tilde{\boldsymbol{\theta}})}{w_k(\tilde{\boldsymbol{\theta}})}\right) \mathbf{M}^{-1}\right)$$

and therefore

(2.8)
$$n^{1/2}((\hat{\boldsymbol{\theta}}_{0,n}/d^2, \hat{\boldsymbol{\theta}}_{1,n}/d^2, \ldots, \hat{\boldsymbol{\theta}}_{p,n}/d^2, \hat{\boldsymbol{\theta}}_{p+1,n}, \ldots, \hat{\boldsymbol{\theta}}_{p+q,n}) - \boldsymbol{\theta})$$
$$\xrightarrow{\mathcal{D}} \mathbf{N}\left(\mathbf{0}, 4\tau^2 E\left(\frac{w'_k(\tilde{\boldsymbol{\theta}})}{w_k(\tilde{\boldsymbol{\theta}})},\right)^T \left(\frac{w'_k(\tilde{\boldsymbol{\theta}})}{w_k(\tilde{\boldsymbol{\theta}})}\right)\right),$$

where $\hat{\boldsymbol{\theta}}_n = (\hat{\boldsymbol{\theta}}_{0,n}/d^2, \hat{\boldsymbol{\theta}}_{1,n}/d^2, \ldots, \hat{\boldsymbol{\theta}}_{p,n}/d^2, \hat{\boldsymbol{\theta}}_{p+1,n}, \ldots, \hat{\boldsymbol{\theta}}_{p+q,n})$. The limit result in (2.8) means that the only term which depends on $h$ in the limit is $\tau = \tau(\tilde{\varepsilon}/d)$. So the efficiency of the estimator is determined by $\tau$ only.

Let us assume that the innovations $\tilde{\varepsilon}_k$ in (2.6) and (2.7) are standard normal random variables. Using the quasi-maximum likelihood method of Example 2.1 (which is the likelihood method of Example 2.4 in this case), we get that $\tau^2_{\text{quasi}} = 1/2$. If we use the method of Example 2.2, we must rescale since it is assumed that the expected value of the absolute value of the innovations is 1, so the standard normal innovation must be divided by $\sqrt{2/\pi}$. Hence $\tau^2_{\text{exp}} = \pi/2 - 1$. Clearly, $\tau^2_{\text{quasi}} < \tau^2_{\text{exp}}$.

Now we assume that the innovations $\tilde{\varepsilon}_k$ are two-sided exponential random variables. In this case the methods of Examples 2.2 and 2.4 are the same and $\tau^2_{\text{exp}} = 1$. If we use the method of Example 2.1, we need that the second moment is 1, so the innovations must be divided by $\sqrt{2}$. Hence $\tau^2_{\text{quasi}} = 5/4$ in Example 2.1. This means that the variance of the estimators for $\beta_1, \ldots, \beta_q$ ($\beta_1, \ldots, \beta_q$ are invariant for rescaling the innovations) will be 25% more if the quasi-maximum likelihood method is used instead of the likelihood method.

Let $\tilde{\varepsilon}_i$ be independent, identically distributed random variables with density function $f(t) = \{(\vartheta - 1)/2\}(1 + |t|)^{-\vartheta}$, where $\vartheta > 5$. Elementary calcu-lations show that $E|\tilde{\varepsilon}_i| = 1/(\vartheta - 2), E|\tilde{\varepsilon}_i|^2 = 2/((\vartheta - 2)(\vartheta - 3))$ and $E|\tilde{\varepsilon}_i|^4 = 24/((\vartheta - 2)(\vartheta - 3)(\vartheta - 4)(\vartheta - 5))$. If the quasi-maximum likelihood method is used, we use $\varepsilon_i = \tilde{\varepsilon}_i/(E|\tilde{\varepsilon}_i|^2)^{1/2}$ and therefore

$$\tau_{\text{quasi}}^2 = \frac{1}{4}\left\{\frac{6(\vartheta - 2)(\vartheta - 3)}{(\vartheta - 4)(\vartheta - 5)} - 1\right\}.$$

If we use the method of Example 2.2, that is, the two-sided exponential density in the definition of $\hat{L}_n(\mathbf{u})$, we use the innovations $\varepsilon_i = \tilde{\varepsilon}_i/E|\tilde{\varepsilon}_i|$ and we get

$$\tau_{\text{exp}}^2 = \frac{2(\vartheta - 2)}{\vartheta - 3} - 1.$$

Elementary calculations show that $\tau_{\text{quasi}}^2 > \tau_{\text{exp}}^2$ for any $\vartheta > 5$. If $\vartheta = 6$, then $\tau_{\text{quasi}}^2 = 8.75$ while $\tau_{\text{exp}}^2 \approx 1.67$. If $\vartheta$ is large then $\tau_{\text{quasi}}^2 \approx 1.25$ while $\tau_{\text{exp}}^2 \approx 1$. The parameters $\beta_1, \ldots, \beta_q$ are invariant for scaling, so the two-sided exponential method gives smaller asymptotic variance than the quasi-maximum method. We note that the likelihood method of Example 2.4 provides the smallest possible variance for the estimation of $\beta_1, \ldots, \beta_q$. However, this example illustrates that if the density is unknown and we suspect that the tail of the distribution of the innovations is polynomial, the two-sided exponential method performs better than the quasi-maximum likelihood.

If we are interested in the estimation of $d$ in the examples above, we can use the residuals. The residuals are defined as $\hat{\varepsilon}_i = y_i/\hat{w}_i(\hat{\boldsymbol{\theta}}_n), 1 < i \leq n$. Let us assume that the estimation is done under the scaling assumption (1.16). Then $\hat{d}_n = (\sum_{1 < i \leq n} \hat{\varepsilon}_i^2/(n-1))^{1/2}$ can be used when we move to a model with scaling assumption $\varepsilon_0^2 = 1$. However, replacing $d$ with $\hat{d}$ in (2.8) will change the asymptotic variance. Using different $h$'s in $\hat{L}_n(\mathbf{u})$, we study models based on different scaling assumptions. Since the parameters in (1.1) and (1.2) are not uniquely defined, scaling assumptions or reparametrizations [cf. Drost and Klaassen (1997) and Newey and Steigerwald (1997)] are required.

**3. Preliminary results.** The first six lemmas are from Berkes, Horváth and Kokoszka (2003).

LEMMA 3.1. *If the conditions of Theorem* 1.1 *are satisfied and* $\sigma_0^2 = w_0(\mathbf{u}^*)$ *with some* $\mathbf{u}^* \in U$, *then* $\mathbf{u}^* = \boldsymbol{\theta}$.

PROOF. This result is part of the proof of Lemma 5.5 in Berkes, Horváth and Kokoszka (2003). □

Let $\log^+ x = \log x$ if $x > 1$ and 0 otherwise.

LEMMA 3.2. *Let $\varphi_0, \varphi_1, \varphi_2, \ldots$ be identically distributed random variables satisfying $E \log^+ |\varphi_0| < \infty$. Then $\sum_{1 \leq k < \infty} \varphi_k z^k$ converges a.s. for all $|z| < 1$.*

LEMMA 3.3. *If (1.5), (1.6) and (1.12) hold, then there is $\delta > 0$ such that*

$$(3.1) \qquad E|y_0^2|^\delta < \infty \quad and \quad E|\sigma_0^2|^\delta < \infty.$$

LEMMA 3.4. *If (1.5), (1.6) and (1.8) hold, then there are constants $0 < C^*, C^{**} < \infty$ and $0 < \rho < 1$ such that*

$$(3.2) \qquad C^* \leq w_k(\mathbf{u}) \leq C^{**}\left(1 + \sum_{1 \leq i < \infty} \rho^i y_{k-i}^2\right), \qquad \mathbf{u} \in U,$$

*and*

$$C^* \leq \hat{w}_k(\mathbf{u}) \leq C^{**}\left(1 + \sum_{1 \leq i < \infty} \rho^i y_{k-i}^2\right), \qquad \mathbf{u} \in U,$$

*for any $-\infty < k < \infty$.*

LEMMA 3.5. *If (1.2), (1.5), (1.6), (1.8), (1.11) and (1.12) hold, then, for any $0 < \kappa^* < \kappa$,*

$$E\left(\sup_{\mathbf{u} \in U} \frac{\sigma_k^2}{w_k(\mathbf{u})}\right)^{\kappa^*} < \infty.$$

LEMMA 3.6. *If (1.5), (1.6), (1.8), (1.11) and (1.12) hold, then*

$$E\left(\sup_{\mathbf{u} \in U} \frac{|w_k'(\mathbf{u})|}{w_k(\mathbf{u})}\right)^{\kappa^*} < \infty,$$

$$E\left(\sup_{\mathbf{u} \in U} \frac{|w_k''(\mathbf{u})|}{w_k(\mathbf{u})}\right)^{\kappa^*} < \infty$$

*and*

$$E \sup_{\mathbf{u} \in U} \left|\frac{w_k'''(\mathbf{u})}{w_k(\mathbf{u})}\right|^{\kappa^*} < \infty$$

*for any $\kappa^* > 0$.*

For any $\mathbf{u} = (x, s_1, \ldots, s_p, t_1, \ldots, t_q) \in U$ and $\gamma > 1$, we define

$$(3.3) \quad U(\gamma, \mathbf{u}) = \left\{\mathbf{u}^* = (x^*, s_1^*, \ldots, s_p^*, t_1^*, \ldots, t_q^*) \in U : \max_{1 \leq j \leq q} t_j^*/t_j \leq \gamma\right\}.$$

LEMMA 3.7. *If (1.5), (1.6), (1.8), (1.11) and (1.12) hold, then for any $-\infty < \kappa^* < \infty$ there is $\gamma > 1$ such that*

$$E\left(\sup\left\{\frac{w_k(\mathbf{u}^*)}{w_k(\mathbf{u})} : \mathbf{u}^* \in U(\gamma', \mathbf{u})\right\}\right)^{\kappa^*} < \infty$$

*for all $\mathbf{u} \in U$ and $1 \leq \gamma' \leq \gamma$.*

PROOF. Due to symmetry we can assume that $\kappa^* > 0$. We note that $\underline{u}/(1-\rho_0) \leq c_0(\mathbf{u})$ for all $\mathbf{u} \in U$ and $0 \leq c_i(\mathbf{u}^*) \leq K_1 \gamma^i c_i(\mathbf{u})$ with some constant $K_1$ for all $\mathbf{u}^* \in U(\gamma, \mathbf{u})$ by Lemma 3.1 of Berkes, Horváth and Kokoszka (2003). Thus Lemma 3.7 will be proven if we show that

$$E\left(\frac{\sum_{1 \leq i < \infty} \gamma^i c_i(\mathbf{u}) y_{k-i}^2}{1 + \sum_{1 \leq i < \infty} c_i(\mathbf{u}) y_{k-i}^2}\right)^{\kappa^*} \leq K_2.$$

By Lemma 3.1 in Berkes, Horváth and Kokoszka (2003) there are constants $c$ and $0 < \rho < 1$ such that

(3.4) $\qquad |c_k(\mathbf{u})| \leq c\rho^k \qquad$ for all $\mathbf{u} \in U$ and $k$.

For any $M \geq 1$ we have

$$\frac{\sum_{1 \leq i < \infty} \gamma^i c_i(\mathbf{u}) y_{k-i}^2}{1 + \sum_{1 \leq i < \infty} c_i(\mathbf{u}) y_{k-i}^2} \leq \gamma^M + \sum_{M < i < \infty} \gamma^i c_i(\mathbf{u}) y_{k-i}^2$$

$$\leq \gamma^M + K_3 \sum_{M < i < \infty} (\gamma\rho)^i y_{k-i}^2,$$

with some $K_3$ on account of (3.4). By the Markov inequality we have

$$P\left\{\sum_{M < i < \infty} (\gamma\rho)^i y_{k-i}^2 > t/2\right\}$$

$$\leq \sum_{M < i < \infty} P\{y_{k-i}^2 > (t/2)(\gamma\rho)^{-i}(1 - (\gamma\rho)^{1/2})(\gamma\rho)^{i/2}\}$$

$$= \sum_{M < i < \infty} P\{|y_0^2|^\delta > (t/2)^\delta (1 - (\gamma\rho)^{1/2})^\delta (\gamma\rho)^{-i\delta/2}\}$$

$$\leq E|y_0^2|^\delta (1 - (\gamma\rho)^{1/2})^{-\delta}(1 - (\gamma\rho)^{\delta/2})^{-1}(t/2)^{-\delta}(\gamma\rho)^{M\delta/2}.$$

Choosing $M = \log(t/2)/\log \gamma$, $t > \gamma^2$, we have, for any $\kappa^* > 0$,

$$P\left\{\sup_{\mathbf{u} \in U} \frac{\sum_{1 \leq i < \infty} \gamma^i c_i(\mathbf{u}) y_{k-i}^2}{1 + \sum_{1 \leq i < \infty} c_i(\mathbf{u}) y_{k-i}^2} > t\right\}$$

$$\leq P\left\{\sum_{M < i < \infty} (\gamma\rho)^i y_{k-i}^2 > t/2\right\}$$

$$\leq K_4 \exp(-(\delta/2)(1+\log\rho^{-1}/\log\gamma)\log(t/2))$$
$$\leq K_5 t^{-2\kappa^*}$$

if $\gamma > 1$ is close enough to 1, where $K_4$ and $K_5$ are constants. This completes the proof of Lemma 3.7. $\square$

## 4. Proofs.

LEMMA 4.1. *If* (1.5), (1.6) *and* (1.8)–(1.13) *hold, then* $L(\mathbf{u})$ *is defined for any* $\mathbf{u} \in U$.

PROOF. By Lemma 3.3 there is $0 < \delta \leq 1$ such that $E|y_0^2|^\delta < \infty$ and therefore

$$
\begin{aligned}
E\left(1 + \sum_{1\leq i<\infty} \rho^i y_{k-i}^2\right)^\delta \\
(4.1) \qquad \leq E\left(1 + \sum_{1\leq i<\infty} \rho^{\delta i} |y_0^2|^\delta\right) \\
= 1 + E|y_0^2|^\delta \sum_{1\leq i<\infty} (\rho^\delta)^i < \infty
\end{aligned}
$$

for all $0 < \rho < 1$. Therefore by Lemma 3.4 we have

$$(4.2) \qquad E\sup_{\mathbf{u}\in U} |\log w_0(\mathbf{u})| < \infty.$$

Since $\varepsilon_0$ and $\sigma_0/w_0^{1/2}(\mathbf{u})$ are independent, by (1.13) we obtain

$$E\left|\log h\left(\varepsilon_0 \frac{\sigma_0}{w_0^{1/2}(\mathbf{u})}\right)\right| \leq C_0\left(1 + E\left(\frac{\sigma_0^2}{w_0(\mathbf{u})}\right)^{\nu_0/2}\right).$$

Using Lemma 3.5 and (4.2), we conclude

$$(4.3) \qquad E\sup_{\mathbf{u}\in U}\left|\log\left\{\frac{1}{w_0^{1/2}(\mathbf{u})}h\left(\frac{y_0}{w_0^{1/2}(\mathbf{u})}\right)\right\}\right| < \infty,$$

and thus Lemma 4.1 is proved. $\square$

Let

$$L_n(\mathbf{u}) = \frac{1}{n}\sum_{1\leq k\leq n} \log\left\{\frac{1}{w_k^{1/2}(\mathbf{u})}h\left(\frac{y_k}{w_k^{1/2}(\mathbf{u})}\right)\right\}.$$

LEMMA 4.2. *If* (1.5), (1.6) *and* (1.8)–(1.14) *hold, then*
$$E\left\{\sup \frac{1}{|\mathbf{u}-\mathbf{v}|}|L_n(\mathbf{u}) - L_n(\mathbf{v})| : \mathbf{u} \in U, \mathbf{v} \in U\right\} < \infty.$$

PROOF. By the mean value theorem there is a random variable $\boldsymbol{\eta} \in U$ such that
$$|g(\varepsilon_k, \sigma_k/w_k^{1/2}(\mathbf{u})) - g(\varepsilon_k, \sigma_k/w_k^{1/2}(\mathbf{v}))|$$
$$= \frac{1}{2}|\mathbf{u}-\mathbf{v}|\left|g_1(\varepsilon_k, \sigma_k/w_k^{1/2}(\boldsymbol{\eta}))\frac{\sigma_k}{w_k^{3/2}(\boldsymbol{\eta})}w_k'(\boldsymbol{\eta})\right|.$$

So by condition (1.14) we conclude
$$\sup_{\mathbf{u},\mathbf{v}} \frac{1}{|\mathbf{u}-\mathbf{v}|}\left|g\left(\varepsilon_k, \frac{\sigma_k}{w_k^{1/2}(\mathbf{u})}\right) - g\left(\varepsilon_k, \frac{\sigma_k}{w_k^{1/2}(\mathbf{v})}\right)\right|$$
$$\leq \frac{1}{2}C_1(\varepsilon_k)\sup_{\boldsymbol{\eta}}\left|\left\{\left(\left(\frac{\sigma_k}{w_k^{1/2}(\boldsymbol{\eta})}\right)^{\nu_1}+1\right)\Big/\left(\frac{\sigma_k}{w_k^{1/2}(\boldsymbol{\eta})}\right)\right\}\frac{\sigma_k}{w_k^{1/2}(\boldsymbol{\eta})}\frac{w_k'(\boldsymbol{\eta})}{w_k(\boldsymbol{\eta})}\right|$$
$$\leq \frac{1}{2}C_1(\varepsilon_k)\left\{\sup_{\mathbf{u}}\left(\left(\frac{\sigma_k^2}{w_k(\mathbf{u})}\right)^{\nu_1/2}+1\right)\right\}\sup_{\mathbf{v}}\left|\frac{w_k'(\mathbf{v})}{w_k(\mathbf{v})}\right|.$$

We note that $\varepsilon_k$ and $\{\sigma_k^2/w_k(\mathbf{u}), \mathbf{u} \in U\}$ are independent for any $k$. The Hölder inequality and Lemmas 3.5 and 3.6 yield
$$EC_1(\varepsilon_k)\left\{\sup_{\mathbf{u}}\left(\left(\frac{\sigma_k^2}{w_k(\mathbf{u})}\right)^{\nu_1/2}+1\right)\right\}\sup_{\mathbf{v}}\left|\frac{w_k'(\mathbf{v})}{w_k(\mathbf{v})}\right|$$
$$\leq EC_1(\varepsilon_k)\left\{E\left(\sup_{\mathbf{u}}\left(\frac{\sigma_k^2}{w_k(\mathbf{u})}\right)^{\nu_1/2}+1\right)^\gamma\right\}^{1/\gamma}\left\{\sup_{\mathbf{v}}\left|\frac{w_k'(\mathbf{v})}{w_k(\mathbf{v})}\right|^{\gamma'}\right\}^{1/\gamma'}$$
$$\leq K_1,$$

with some constant $K_1$, where $1 < \gamma, \gamma' < \infty$ satisfy $(\nu_1/2)\gamma < \kappa$ and $1/\gamma + 1/\gamma' = 1$. Since $L_n(\mathbf{u})$ is the average of stationary random variables, the proof of Lemma 4.2 is now complete. $\square$

LEMMA 4.3. *If* (1.5), (1.6) *and* (1.8)–(1.14) *hold, then*
$$n\sup_{\mathbf{u}\in U}|\hat{L}_n(\mathbf{u}) - L_n(\mathbf{u})| = O(1) \qquad a.s.$$

PROOF. We use (3.4). Let
$$\xi = \sum_{1\leq i<\infty}\rho^i y_{-i}^2.$$

We note that by Lemmas 3.2 and 3.3 the series defining $\xi$ converges a.s. Using the definitions of $w_k(\mathbf{u})$, $\hat{w}_k(\mathbf{u})$ and (3.4), we conclude
(4.4)
$$\sup_{\mathbf{u}\in U}|w_k(\mathbf{u}) - \hat{w}_k(\mathbf{u})| \leq c\sum_{k<i<\infty}\rho^i y_{k-i}^2 = c\rho^k \xi.$$

Using next the mean value theorem, there is $\eta \in (\sigma_k/\hat{w}_k^{1/2}(\mathbf{u}), \sigma_k/w_k^{1/2}(\mathbf{u}))$ such that

$$|g(\varepsilon_k, \sigma_k/w_k^{1/2}(\mathbf{u})) - g(\varepsilon_k, \sigma_k/\hat{w}_k^{1/2}(\mathbf{u}))|$$
$$= \left| g_1(\varepsilon_k, \eta) \left( \frac{\sigma_k}{w_k^{1/2}(\mathbf{u})} - \frac{\sigma_k}{\hat{w}_k^{1/2}(\mathbf{u})} \right) \right|.$$

Applying condition (1.14) and Lemma 3.4, we get

$$\left| g_1(\varepsilon_k, \eta) \left( \frac{\sigma_k}{w_k^{1/2}(\mathbf{u})} - \frac{\sigma_k}{\hat{w}_k^{1/2}(\mathbf{u})} \right) \right|$$
$$\leq C_1(\varepsilon_k) \{ (\eta^{\nu_1} + 1)/\eta \} \sigma_k \left| \frac{w_k(\mathbf{u}) - \hat{w}_k(\mathbf{u})}{w_k^{1/2}(\mathbf{u}) \hat{w}_k^{1/2}(\mathbf{u}) (w_k^{1/2}(\mathbf{u}) + \hat{w}_k^{1/2}(\mathbf{u}))} \right|$$
$$\leq C_1(\varepsilon_k) \left\{ \left( \frac{\sigma_k}{\hat{w}_k^{1/2}(\mathbf{u})} \right)^{\nu_1} + 1 \right\} \frac{w_k^{1/2}(\mathbf{u})}{\sigma_k} \cdot \frac{\sigma_k}{w_k^{1/2}(\mathbf{u})} \cdot \frac{|w_k(\mathbf{u}) - \hat{w}_k(\mathbf{u})|}{2C^*}$$
$$\leq K_2 C_1(\varepsilon_k)(\sigma_k^{\nu_1} + 1) \xi \rho^k$$

for any $k$. Applying (4.2) with $\mathbf{u} = \boldsymbol{\theta}$, we see that $E|\log \sigma_0| < \infty$. We can assume without loss of generality that $C_1$ in (1.14) is larger than 1 and thus by (1.15) we conclude that $E|\log(C_1(\varepsilon_0)(\sigma_0^{\nu_1} + 1))| < \infty$. Thus we can apply Lemma 3.2 to get

$$n \sup_{\mathbf{u} \in U} |L_n(\mathbf{u}) - \hat{L}_n(\mathbf{u})| \leq K_2 \xi \sum_{1 \leq k < \infty} C_1(\varepsilon_k)(\sigma_k^{\nu_1} + 1)\rho^k < \infty \quad \text{a.s.} \qquad \square$$

PROOF OF THEOREM 1.1. Since $\log(w_k^{-1/2}(\mathbf{u})h(y_k w_k^{-1/2}(\mathbf{u})))$ is a stationary sequence with finite mean $L(\mathbf{u})$ and by Theorem 3.5.8 of Stout (1974) it is also ergodic, the ergodic theorem implies that $L_n(\mathbf{u}) \to L(\mathbf{u})$ a.s. for any fixed $\mathbf{u} \in U$. Thus Lemma 4.2 yields

$$\sup_{\mathbf{u} \in U} |L_n(\mathbf{u}) - L(\mathbf{u})| \to 0 \qquad \text{a.s.}$$

Using now Lemma 4.3, we conclude that

(4.5) $$\sup_{\mathbf{u} \in U} |\hat{L}_n(\mathbf{u}) - L(\mathbf{u})| \to 0 \qquad \text{a.s.}$$

We note that

(4.6) $$L(\boldsymbol{\theta}) - L(\mathbf{u}) = E\{g(\varepsilon_0, 1) - g(\varepsilon_0, \sigma_0/w^{1/2}(\mathbf{u}))\}.$$

Since $\varepsilon_0$ and $\{\sigma_0/w_0^{1/2}(\mathbf{u}), \mathbf{u} \in U\}$ are independent, by (1.16) we have that $L(\boldsymbol{\theta}) \geq L(\mathbf{u})$ and we have that $L(\boldsymbol{\theta}) = L(\mathbf{u})$ if and only if $\sigma_0 = w_0^{1/2}(\mathbf{u})$. Using Lemma 3.1, we get that $L(\mathbf{u})$ has a unique maximum at $\boldsymbol{\theta}$. The function $L(\mathbf{u})$ is continuous, and thus the uniform a.s. convergence of $\hat{L}_n(\mathbf{u})$ to $L(\mathbf{u})$ implies $\hat{\boldsymbol{\theta}}_n \to \boldsymbol{\theta}$, proving Theorem 1.1.

Since
$$L_n(\mathbf{u}) = \frac{1}{n} \sum_{1 \leq k \leq n} g(\varepsilon_k, \sigma_k/w_k^{1/2}(\mathbf{u})) - \frac{1}{n} \sum_{1 \leq k \leq n} \log \sigma_k,$$

we get that

$$(4.7) \quad L'_n(\mathbf{u}) = \frac{1}{n} \sum_{1 \leq k \leq n} g_1(\varepsilon_k, \sigma_k/w_k^{1/2}(\mathbf{u})) \left( -\frac{1}{2} \frac{\sigma_k}{w_k^{1/2}(\mathbf{u})} \frac{w'_k(\mathbf{u})}{w_k(\mathbf{u})} \right)$$

and

$$L''_n(\mathbf{u}) = \frac{1}{n} \sum_{1 \leq k \leq n} g_2(\varepsilon_k, \sigma_k/w_k^{1/2}(\mathbf{u}))$$
$$\times \left( -\frac{1}{2} \frac{\sigma_k}{w_k^{1/2}(\mathbf{u})} \frac{w'_k(\mathbf{u})}{w_k(\mathbf{u})} \right)^T \left( -\frac{1}{2} \frac{\sigma_k}{w_k^{1/2}(\mathbf{u})} \frac{w'_k(\mathbf{u})}{w_k(\mathbf{u})} \right)$$
$$(4.8) \quad + \frac{1}{n} \sum_{1 \leq k \leq n} g_1(\varepsilon_k, \sigma_k/w_k^{1/2}(\mathbf{u}))$$
$$\times \left( \frac{3}{4} \frac{\sigma_k}{w_k^{1/2}(\mathbf{u})} \left( \frac{w'_k(\mathbf{u})}{w_k(\mathbf{u})} \right)^T \frac{w'_k(\mathbf{u})}{w_k(\mathbf{u})} - \frac{1}{2} \frac{\sigma_k}{w_k^{1/2}(\mathbf{u})} \frac{w''_k(\mathbf{u})}{w_k(\mathbf{u})} \right).$$

Similarly,

$$(4.9) \quad L'(\mathbf{u}) = E g_1(\varepsilon_0, \sigma_0/w_0^{1/2}(\mathbf{u})) \left( -\frac{1}{2} \frac{\sigma_0}{w_0^{1/2}(\mathbf{u})} \frac{w'_0(\mathbf{u})}{w_0(\mathbf{u})} \right)$$

and

$$L''(\mathbf{u}) = E \Bigg\{ g_2(\varepsilon_0, \sigma_0/w_0^{1/2}(\mathbf{u})) \left( -\frac{1}{2} \frac{\sigma_0}{w_0^{1/2}(\mathbf{u})} \frac{w'_0(\mathbf{u})}{w_0(\mathbf{u})} \right)^T \left( -\frac{1}{2} \frac{\sigma_0}{w_0^{1/2}(\mathbf{u})} \frac{w'_0(\mathbf{u})}{w_0(\mathbf{u})} \right)$$
$$(4.10) \quad + g_1(\varepsilon_0, \sigma_0/w_0^{1/2}(\mathbf{u})) \left( \frac{3}{4} \frac{\sigma_0}{w_0^{1/2}(\mathbf{u})} \left( \frac{w'_0(\mathbf{u})}{w_0(\mathbf{u})} \right)^T \frac{w'_0(\mathbf{u})}{w_0(\mathbf{u})} \right.$$
$$\left. - \frac{1}{2} \frac{\sigma_0}{w_0^{1/2}(\mathbf{u})} \frac{w''_0(\mathbf{u})}{w_0(\mathbf{u})} \right) \Bigg\}.$$

The expected value in (4.9) exists, since by (1.14) and (1.15) and the independence of $\varepsilon_0$ and $\sigma_0/w^{1/2}(\mathbf{u})$ we have

$$E \left| g_1 \left( \varepsilon_0, \frac{\sigma_0}{w_0^{1/2}(\mathbf{u})} \right) \left( \frac{\sigma_0}{w_0^{1/2}(\mathbf{u})} \frac{w'_0(\mathbf{u})}{w_0(\mathbf{u})} \right) \right|$$
$$(4.11) \quad \leq E C_1(\varepsilon_0) \left( \left( \frac{\sigma_0}{w_0^{1/2}(\mathbf{u})} \right)^{\nu_1} + 1 \right) \left| \frac{w'_0(\mathbf{u})}{w_0(\mathbf{u})} \right|$$
$$\leq E C_1(\varepsilon_0) E \left( \left( \frac{\sigma_0}{w_0^{1/2}(\mathbf{u})} \right)^{\nu_1} + 1 \right) \left| \frac{w'_0(\mathbf{u})}{w_0(\mathbf{u})} \right| < \infty$$

on account of the Hölder inequality and Lemmas 3.5 and 3.6. A similar argument shows that the expected value in (4.10) also exists for all $\mathbf{u} \in U$. □

LEMMA 4.4. *If* (1.5), (1.6), (1.8)–(1.14), (1.17) *and* (1.18) *hold, then there exists* $U^*$, *a neighborhood of* $\boldsymbol{\theta}$, *such that*

$$
\sup_{\mathbf{u} \in U^*} |L_n'(\mathbf{u}) - L'(\mathbf{u})| \to 0 \qquad a.s. \tag{4.12}
$$

*If, in addition,* (1.19) *and* (1.20) *are satisfied, then*

$$
\sup_{\mathbf{u} \in U^*} |L_n''(\mathbf{u}) - L''(\mathbf{u})| \to 0 \qquad a.s. \tag{4.13}
$$

*Also,* $E\{(w_0'(\mathbf{u})/w_0(\mathbf{u}))^T w_0'(\mathbf{u})/w_0(\mathbf{u})\}$ *is a nonsingular matrix for any* $\mathbf{u} \in U^*$.

PROOF. Let $U^* = U(\gamma, \boldsymbol{\theta})$ with some $\gamma$ to be chosen later. Applying (1.17), we obtain

$$
\begin{aligned}
\left| \left( g_1\left(\varepsilon_k, \frac{\sigma_k}{w_k^{1/2}(\mathbf{u})}\right) \right)' \right| \\
= \left| g_2\left(\varepsilon_k, \frac{\sigma_k}{w_k^{1/2}(\mathbf{u})}\right) \sigma_k \left(-\frac{1}{2}\right) \frac{w_k'(\mathbf{u})}{w_k^{3/2}(\mathbf{u})} \right| \\
\leq C_2(\varepsilon_k) \left( \left( \frac{\sigma_k}{w_k^{1/2}(\mathbf{u})} \right)^{\nu_2} + 1 \right) \left( \frac{\sigma_k}{w_k^{1/2}(\mathbf{u})} \right)^{-1} \left| \frac{w_k'(\mathbf{u})}{w_k(\mathbf{u})} \right|.
\end{aligned}
\tag{4.14}
$$

Using (4.7), (4.14) and conditions (1.14) and (1.17), we get

$$
\frac{1}{|\boldsymbol{\theta} - \mathbf{u}|} n |L_n'(\boldsymbol{\theta}) - L_n'(\mathbf{u})|
$$

$$
\leq \frac{1}{|\boldsymbol{\theta} - \mathbf{u}|} \sum_{1 \leq k \leq n} \left| g_1\left(\varepsilon_k, \frac{\sigma_k}{w_k^{1/2}(\boldsymbol{\theta})}\right) - g_1\left(\varepsilon_k, \frac{\sigma_k}{w_k^{1/2}(\mathbf{u})}\right) \right| \left| \frac{\sigma_k}{w_k^{1/2}(\boldsymbol{\theta})} \frac{w_k'(\boldsymbol{\theta})}{w_k(\boldsymbol{\theta})} \right|
$$

$$
+ \frac{1}{|\boldsymbol{\theta} - \mathbf{u}|} \sum_{1 \leq k \leq n} \left| g_1\left(\varepsilon_k, \frac{\sigma_k}{w_k^{1/2}(\mathbf{u})}\right) \right| \left| \frac{\sigma_k}{w_k^{1/2}(\boldsymbol{\theta})} \frac{w_k'(\boldsymbol{\theta})}{w_k(\boldsymbol{\theta})} - \frac{\sigma_k}{w_k^{1/2}(\mathbf{u})} \frac{w_k'(\mathbf{u})}{w_k(\mathbf{u})} \right|
$$

$$
\leq \sum_{1 \leq k \leq n} C_2(\varepsilon_k) \left( \sup_{\mathbf{v} \in U^*} \left( \left( \frac{\sigma_k}{w_k^{1/2}(\mathbf{v})} \right)^{\nu_2} + 1 \right) \left( \frac{\sigma_k}{w_k^{1/2}(\mathbf{v})} \right)^{-1} \right) \sup_{\mathbf{u} \in U} \left| \frac{w_k'(\mathbf{u})}{w_k(\mathbf{u})} \right|^2
$$

$$
+ \sum_{1 \leq k \leq n} C_1(\varepsilon_k) \left( \sup_{\mathbf{v} \in U^*} \left( \left( \frac{\sigma_k}{w_k^{1/2}(\mathbf{v})} \right)^{\nu_1} + 1 \right) \left( \frac{\sigma_k}{w_k^{1/2}(\mathbf{v})} \right) \right)^{-1}
$$

$$\times \sup_{\mathbf{z}\in U^*} \frac{\sigma_k}{w_k^{1/2}(\mathbf{z})} \sup_{\mathbf{u}\in U} \left(\frac{3}{2}\left|\frac{w_k'(\mathbf{u})}{w_k(\mathbf{u})}\right|^2 + \left|\frac{w_k''(\mathbf{u})}{w_k(\mathbf{u})}\right|\right)$$

$$= \sum_{1\le k\le n}(I_{k,1}+I_{k,2}).$$

Using the Cauchy inequality and the independence of $\varepsilon_k$ and $\{w_k(\mathbf{u}),\mathbf{u}\in U\}$, we conclude from Lemmas 3.7 and 3.6 that

$$EI_{k,1} = EC_2(\varepsilon_k)E\left(\sup_{\mathbf{v}\in U^*}\left(\left(\frac{\sigma_k}{w_k^{1/2}(\mathbf{v})}\right)^{\nu_2}+1\right)\left(\frac{\sigma_k}{w_k^{1/2}(\mathbf{v})}\right)^{-1}\sup_{\mathbf{u}\in U}\left|\frac{w_k'(\mathbf{u})}{w_k(\mathbf{u})}\right|^2\right)$$

$$\le EC_2(\varepsilon_k)\left(E\left(\sup_{\mathbf{v}\in U^*}\left(\left(\frac{\sigma_k}{w_k^{1/2}(\mathbf{v})}\right)^{\nu_2}+1\right)\left(\frac{\sigma_k}{w_k^{1/2}(\mathbf{v})}\right)^{-1}\right)^2\right)^{1/2}$$

$$\times \left(E\sup_{\mathbf{u}\in U}\left|\frac{w_k'(\mathbf{u})}{w_k(\mathbf{u})}\right|^4\right)^{1/2} < \infty,$$

provided $\gamma > 1$ is chosen close enough to 1. Similarly,

$$EI_{k,2} < \infty.$$

Thus $I_{k,1}$ and $I_{k,2}$ are stationary sequences with finite expectations and by Theorem 3.4.8 of Stout (1974) they are also ergodic. Hence our previous estimates and the ergodic theorem imply

(4.15) $$\sup_{\mathbf{u}\in U^*}\frac{1}{|\boldsymbol{\theta}-\mathbf{u}|}|L_n'(\boldsymbol{\theta})-L_n'(\mathbf{u})| = O(1) \qquad \text{a.s.}$$

Since $L_n'(\mathbf{u})$ is an average of a stationary, ergodic sequence with finite expectation [cf. (4.11)], another application of the ergodic theorem gives, for any $\mathbf{u}\in U^*$,

(4.16) $$L_n'(\mathbf{u}) \to L'(\mathbf{u}) \qquad \text{a.s.}$$

Putting together (4.15) and (4.16), we get (4.12). Similar arguments yield (4.13).

Berkes, Horváth and Kokoszka (2003) proved that $E(w_0'(\mathbf{u}))^T w_0'(\mathbf{u})/w_0^2(\mathbf{u})$ is a continuous function and it is nonsingular at $\mathbf{u}=\boldsymbol{\theta}$, so the proof of Lemma 4.4 is complete. $\square$

Let $\boldsymbol{\theta}_n = \arg\max\{L_n(\mathbf{u}):\mathbf{u}\in U\}$.

LEMMA 4.5. *If* (1.5), (1.6) *and* (1.8)–(1.22) *are satisfied, then*

(4.17)
$$n^{1/2}(\boldsymbol{\theta}_n - \boldsymbol{\theta})$$
$$= \frac{1}{n^{1/2}}\sum_{1\le k\le n} g_1(\varepsilon_k,1)\frac{2}{Eg_2(\varepsilon_0,1)}\frac{w_k'(\boldsymbol{\theta})}{w_k(\boldsymbol{\theta})}\mathbf{A}^{-1}(1+o(1)) \qquad a.s.$$

*and*

$$n^{1/2}(\boldsymbol{\theta}_n - \boldsymbol{\theta})$$

(4.18)
$$= \frac{1}{n^{1/2}} \sum_{1 \leq k \leq n} g_1(\varepsilon_k, 1) \frac{2}{Eg_2(\varepsilon_0, 1)} \frac{w'_k(\boldsymbol{\theta})}{w_k(\boldsymbol{\theta})} \mathbf{A}^{-1} + o_P(1).$$

PROOF. We showed [cf. (4.9) and (4.11)] that $L(\mathbf{u})$ is differentiable for all $\mathbf{u} \in U$ and proved after (4.6) that $L(\mathbf{u})$ has a unique maximum at $\mathbf{u} = \boldsymbol{\theta}$. Thus $L'(\boldsymbol{\theta}) = \mathbf{0}$. By the independence of $\varepsilon_0$ and $w'_0(\boldsymbol{\theta})/w_0(\boldsymbol{\theta})$ and (4.9), we have

$$L'(\boldsymbol{\theta}) = Eg_1(\varepsilon_0, 1)\left(-\frac{1}{2} E \frac{w'_0(\boldsymbol{\theta})}{w_0(\boldsymbol{\theta})}\right).$$

Berkes, Horváth and Kokoszka (2003) showed $Ew'_0(\boldsymbol{\theta})/w_0(\boldsymbol{\theta}) \neq \mathbf{0}$, so we have

(4.19)
$$Eg_1(\varepsilon_0, 1) = 0.$$

By (4.5) it follows easily that $\boldsymbol{\theta}_n \to \boldsymbol{\theta}$ a.s. Hence there is a random variable $n_0$ such that $\boldsymbol{\theta}_n \in U^*$ if $n \geq n_0$, where $U^*$ is defined in Lemma 4.4. Clearly, $U^* \subset U$ is compact and for $\gamma > 1$ sufficiently close to 1 it does not have common points with the boundary of $U$. Since $L_n(\mathbf{u})$ is twice differentiable and it reaches a maximum at $\boldsymbol{\theta}_n$, we have

(4.20)
$$L'_n(\boldsymbol{\theta}_n) = \mathbf{0} \qquad \text{if } n \geq n_0,$$

and thus

$$L'_n(\boldsymbol{\theta}_n) - L'_n(\boldsymbol{\theta}) = -L'_n(\boldsymbol{\theta}).$$

By (4.13) we have that $L'_n(\boldsymbol{\theta}_n) - L'_n(\boldsymbol{\theta}) = (\boldsymbol{\theta}_n - \boldsymbol{\theta})(L''(\boldsymbol{\theta}) + o(1))$ a.s. Observing that $L''(\boldsymbol{\theta}) = Eg_2(\varepsilon_0, 1)\frac{1}{4}\mathbf{A}$, and using (4.7), the proof of (4.17) is complete. By the orthogonality and stationarity of the summands in (4.17), and in view of (1.21), the variance of the sum is $O(1/n)$ and therefore (4.18) follows from (4.17). □

A simple calculation shows, in analogy with (4.7),

(4.21) $\quad \hat{L}'_n(\mathbf{u}) = \frac{1}{n} \sum_{1 \leq k \leq n} g_1(\varepsilon_k, \sigma_k/\hat{w}_k^{1/2}(\mathbf{u}))\left(-\frac{1}{2} \frac{\sigma_k}{\hat{w}_k^{1/2}(\mathbf{u})} \frac{\hat{w}'_k(\mathbf{u})}{\hat{w}_k(\mathbf{u})}\right).$

LEMMA 4.6. *If* (1.5), (1.6) *and* (1.8)–(1.22) *are satisfied, then*

(4.22)
$$\sup_{\mathbf{u} \in U}|\hat{L}'_n(\mathbf{u}) - L'_n(\mathbf{u})| = O\left(\frac{1}{n}\right) \qquad a.s.$$

PROOF. Berkes, Horváth and Kokoszka (2003) showed that there are constants $c$ and $0 < \rho_* < 1$ such that

$$|c_i'(\mathbf{u})| \leq c\rho_*^i \quad \text{and} \quad |c_i''(\mathbf{u})| \leq c\rho_*^i$$

for all $\mathbf{u} \in U$ and $0 \leq i < \infty$. Hence

(4.23) $$\sup_{\mathbf{u} \in U} |w_k'(\mathbf{u}) - \hat{w}_k'(\mathbf{u})| \leq c \sum_{k < i < \infty} \rho_*^i y_{k-i}^2 = c\rho_*^k \xi_*,$$

where $\xi_* = \sum_{1 \leq i < \infty} \rho_*^i y_{k-i}^2$ converges a.s. by Lemmas 3.2 and 3.3. By (1.14), (1.17), (4.7), (4.21) and $\hat{w}_k(\mathbf{u}) \leq w_k(\mathbf{u})$, we have

$$n|\hat{L}_n'(\mathbf{u}) - L_n(\mathbf{u})|$$

$$\leq \sum_{1 \leq k \leq n} \left|g_1\left(\varepsilon_k, \frac{\sigma_k}{\hat{w}_k^{1/2}(\mathbf{u})}\right) - g_1\left(\varepsilon_k, \frac{\sigma_k}{w_k^{1/2}(\mathbf{u})}\right)\right| \left|\frac{\sigma_k}{w_k^{1/2}(\mathbf{u})} \frac{w_k'(\mathbf{u})}{w_k(\mathbf{u})}\right|$$

$$+ \sum_{1 \leq k \leq n} \left|g_1\left(\varepsilon_k, \frac{\sigma_k}{\hat{w}_k^{1/2}(\mathbf{u})}\right)\right| \left|\frac{\sigma_k}{\hat{w}_k^{1/2}(\mathbf{u})} \frac{\hat{w}_k'(\mathbf{u})}{\hat{w}_k(\mathbf{u})} - \frac{\sigma_k}{w_k^{1/2}(\mathbf{u})} \frac{w_k'(\mathbf{u})}{w_k(\mathbf{u})}\right|$$

$$\leq \sum_{1 \leq k \leq n} C_2(\varepsilon_k) \left(\left(\frac{\sigma_k}{\hat{w}_k^{1/2}(\mathbf{u})}\right)^{\nu_2} + 1\right) \left(\frac{\sigma_k}{w_k^{1/2}(\mathbf{u})}\right)^{-2}$$

$$\times \left|\frac{\sigma_k}{w_k^{1/2}(\mathbf{u})} - \frac{\sigma_k}{\hat{w}_k^{1/2}(\mathbf{u})}\right| \left|\frac{\sigma_k}{w_k^{1/2}(\mathbf{u})} \frac{w_k'(\mathbf{u})}{w_k(\mathbf{u})}\right|$$

$$+ \sum_{1 \leq k \leq n} C_1(\varepsilon_k) \left(\left(\frac{\sigma_k}{\hat{w}_k^{1/2}(\mathbf{u})}\right)^{\nu_1} + 1\right) \left(\frac{\sigma_k}{\hat{w}_k^{1/2}(\mathbf{u})}\right)^{-1}$$

$$\times \sigma_k \left|\frac{\hat{w}_k'(\mathbf{u})}{\hat{w}_k^{3/2}(\mathbf{u})} - \frac{w_k'(\mathbf{u})}{w_k^{3/2}(\mathbf{u})}\right|$$

$$= J_{n,1}(\mathbf{u}) + J_{n,2}(\mathbf{u}).$$

By Lemma 3.4 and (4.4), using $\xi$ and $\rho$ in the proof of Lemma 4.3, we get

$$J_{n,1}(\mathbf{u}) \leq \sum_{1 \leq k \leq n} C_2(\varepsilon_k) \left(\left(\frac{\sigma_k}{\hat{w}_k^{1/2}(\mathbf{u})}\right)^{\nu_2} + 1\right) \left(\frac{\sigma_k}{w_k^{1/2}(\mathbf{u})}\right)^{-1}$$

$$\times \left|\frac{w_k'(\mathbf{u})}{w_k(\mathbf{u})}\right| \left|\frac{\sigma_k}{w_k^{1/2}(\mathbf{u})} - \frac{\sigma_k}{\hat{w}_k^{1/2}(\mathbf{u})}\right|$$

$$\leq \sum_{1 \leq k \leq n} C_2(\varepsilon_k) \left(\left(\frac{\sigma_k}{\hat{w}_k^{1/2}(\mathbf{u})}\right)^{\nu_2} + 1\right)$$

$$\times \left|\frac{w_k'(\mathbf{u})}{w_k(\mathbf{u})}\right| \left|\frac{w_k(\mathbf{u}) - \hat{w}_k(\mathbf{u})}{\hat{w}_k^{1/2}(\mathbf{u})(\hat{w}_k^{1/2}(\mathbf{u}) + w_k^{1/2}(\mathbf{u}))}\right|$$

$$\leq K_3 \sum_{1\leq k\leq n} C_2(\varepsilon_k)(\sigma_k^{\nu_2}+1)\sup_{\mathbf{u}\in U}\left|\frac{w'_k(\mathbf{u})}{w_k(\mathbf{u})}\right|\sup_{\mathbf{u}\in U}|w_k(\mathbf{u})-\hat{w}_k(\mathbf{u})|$$

$$\leq K_3\xi\sum_{1\leq k\leq n} C_2(\varepsilon_k)(\sigma_k^{\nu_2}+1)\sup_{\mathbf{u}\in U}\left|\frac{w'_k(\mathbf{u})}{w_k(\mathbf{u})}\right|\rho^k$$

$$\leq K_3\xi\sum_{1\leq k<\infty} C_2(\varepsilon_k)(\sigma_k^{\nu_2}+1)\sup_{\mathbf{u}\in U}\left|\frac{w'_k(\mathbf{u})}{w_k(\mathbf{u})}\right|\rho^k < \infty \qquad \text{a.s.}$$

In the last step we also used Lemma 3.2 and the observation that $C_2(\varepsilon_k)(\sigma_k^{\nu_2}+1) \times \sup_{\mathbf{u}\in U}|w'_k(\mathbf{u})/w_k(\mathbf{u})|$ is a stationary sequence with

$$E\left|\log\left(C_2(\varepsilon_k)(\sigma_k^{\nu_2}+1)\sup_{\mathbf{u}\in U}\left|\frac{w'_k(\mathbf{u})}{w_k(\mathbf{u})}\right|\right)\right| < \infty.$$

Hence $\sup_{\mathbf{u}\in U} J_{n,1} = O(1)$ a.s. Replacing (4.4) with (4.23), similar arguments show that $\sup_{\mathbf{u}\in U} J_{n,2} = O(1)$ a.s., completing the proof of (4.22). □

LEMMA 4.7. *If* (1.5), (1.6) *and* (1.8)–(1.22) *are satisfied, then*

(4.24) $$|\hat{\boldsymbol{\theta}}_n - \boldsymbol{\theta}_n| = O\left(\frac{1}{n}\right) \qquad a.s.$$

PROOF. Similarly to the proof of Lemma 4.5 there is a random variable $n_0$ such that

(4.25) $$\hat{L}'_n(\hat{\boldsymbol{\theta}}_n) = \mathbf{0} \quad \text{and} \quad \hat{\boldsymbol{\theta}}_n \in U^* \qquad \text{if } n\geq n_0,$$

where the set $U^*$ is defined in the proof of Lemma 4.5. By (4.13) we have

$$L'_n(\boldsymbol{\theta}_n) - L'_n(\hat{\boldsymbol{\theta}}_n) = (\boldsymbol{\theta}_n - \hat{\boldsymbol{\theta}}_n)L''(\boldsymbol{\theta})(1+o(1)) \qquad \text{a.s.}$$

and therefore

$$(\boldsymbol{\theta}_n - \hat{\boldsymbol{\theta}}_n) = (L'_n(\boldsymbol{\theta}_n) - L'_n(\hat{\boldsymbol{\theta}}_n))(L''(\boldsymbol{\theta}))^{-1}(1+o(1)) \qquad \text{a.s.}$$

We recall that $L'_n(\boldsymbol{\theta}_n) = \mathbf{0}$. Lemma 4.6 and (4.25) yield that $L'_n(\hat{\boldsymbol{\theta}}_n) = \hat{L}'_n(\hat{\boldsymbol{\theta}}_n) + O(1/n) = O(1/n)$ a.s., completing the proof of Lemma 4.7. □

PROOF OF THEOREM 1.2. By Lemma 4.7, relation (4.18) remains valid if we replace $\boldsymbol{\theta}_n$ by $\hat{\boldsymbol{\theta}}_n$. Observe now that $g_1(\varepsilon_k,1)w'_k(\boldsymbol{\theta})/w_k(\boldsymbol{\theta})$ is a stationary martingale difference sequence with respect to the $\sigma$-algebra generated by $\{\varepsilon_j, j<k\}$. By Theorem 3.4.8 of Stout (1974) it is also ergodic. Using the Cramér–Wold device [cf. Billingsley (1968), page 49] and Theorem 23.1 of Billingsley [(1968), page 206], we obtain the multivariate central limit theorem expressed by Theorem 1.2. □

A. RÉNYI INSTITUTE OF MATHEMATICS
HUNGARIAN ACADEMY OF SCIENCES
P.O. BOX 127
H-1364 BUDAPEST
HUNGARY
E-MAIL: berkes@renyi.hu

DEPARTMENT OF MATHEMATICS
UNIVERSITY OF UTAH
155 SOUTH 1440 EAST
SALT LAKE CITY, UTAH 84112-0090
USA
E-MAIL: horvath@math.utah.edu